%% file: bem38.tex
\documentclass[royal,cmfonts,final]{witpress}

\input{preamble.tex}

\begin{document}

\title{Multiscale models and approximation algorithms for protein electrostatics}

\author{J. P. Bardhan${}^1$ and M. G. Knepley${}^2$}
\address{${}^1$Dept.~of Mechanical and Industrial Engineering, Northeastern University, Boston MA 02115\\
${}^2$Department of Computational and Applied Mathematics, Rice University, Houston TX 77005}
\maketitle

\section*{Abstract}
Electrostatic forces play many important roles in molecular biology,
but are hard to model due to the complicated interactions between biomolecules and the surrounding solvent, a fluid composed of water and dissolved ions.  Continuum model have been surprisingly successful for simple biological questions, but fail for important problems such as understanding the effects of protein mutations.  In this paper we highlight the advantages of boundary-integral methods for these problems, and our use of boundary integrals to design and test more accurate theories.  Examples include a multiscale model based on nonlocal continuum theory, and a nonlinear boundary condition that captures atomic-scale effects at biomolecular surfaces.

\vspace{1cm}
\textbf{Keywords: electrostatics, proteins, solvation, multiscale, nonlocal, nonlinear, electrolyte, boundary-integral equations, boundary-element methods}

\section{Introduction}\label{sec:introduction}
The behavior of biomolecules such as proteins and DNA depends strongly on their electrostatic interactions with the surrounding fluid, a solvent made of water and dissolved ions~\cite{Bardhan12_review}. Rigorous statistical mechanical theory allows one to use continuum mathematics instead of much slower and more complicated atomistic simulations with explicit solvent molecules~\cite{Roux99}; most continuum theories rely essentially on partial-differential equations (PDEs) based on the Poisson equation and macroscopic dielectric theory~\cite{Bardhan12_review,Roux99}. For numerous modeling problems, e.g. modeling the solute biomolecule using quantum mechanics, boundary-integral equation (BIE) methods enjoy the usual advantages over PDE solvers~\cite{Miertus81}. In this paper, we highlight emerging areas in biomolecular modeling that are enabled by BIEs and fast boundary-element method (BEM) simulation.

The next section introduces common BIE approaches for solving continuum electrostatics.  The following sections then describe our work improving model realism while preserving BEM speed advantages. In particular, we have implemented multiscale models using nonlocal dielectrics (Section~\ref{sec:nonlocal}), and a nonlinear boundary condition model for atomic-scale phenomena at the biomolecule--solvent boundary (Section~\ref{sec:nlbc}).  To encourage discussion and participation by the broader community, each section highlights open questions. Section~\ref{sec:discussion} concludes the paper with a discussion. Space constraints limit our bibliography here, and we welcome interested readers to contact us or consult the more extensive references in our recent reviews~\cite{Bardhan12_review,Bardhan13_nonlocal_review,Bardhan15_CEBA}.

\section{Background}\label{sec:background}

The basic continuum model for understanding protein-solvent electrostatics treats the protein and solvent as distinct media with an interface separating the two regions.  The exterior solvent region (region $I$) is modeled as a continuum dielectric (permittivity $\epsilon_{w} \approx 80$) where the potential obeys either the Laplace equation $\nabla^2 \varphi_{I}(r) = 0$ (modeling pure water) or the linearized Poisson--Boltzmann equation $\nabla^2 \varphi_{I}(r) = \kappa^2 \varphi_{I}(r)$, modeling dilute ionic solution ($\kappa$ is the inverse Debye screening length~\cite{Sharp90}).  The protein (region $II$) is treated as a low-dielectric continuum (relative permittivity $\epsilon_{p} \approx 2-4$) containing an embedded charge distribution $\rho(r)$, where the potential obeys the Poisson equation $\nabla^2 \varphi_{II}(r) = -\rho(r)/\epsilon_0$.  The potential is assumed to decay sufficiently fast as $|r|\rightarrow \infty$, and at the interface $\Gamma$, the permittivity is discontinous, and the potential and normal displacement field are continuous (i.e., $\varphi_{I}(r_\Gamma) = \varphi_{II}(r_\Gamma)$ and $\epsilon_{w}\frac{\varphi_{I}}{\partial n}(r_\Gamma) = \epsilon_{p}\frac{\varphi_{II}}{\partial n}(r_\Gamma)$).  Various BIE formulations can be written (for history and analysis, see~\cite{Bardhan09_disc,Bardhan12_review}).  When the Laplace equation governs in the solvent, one may use the second-kind BIE
\begin{align}
  \left(I  + \hat{\epsilon} \left(-\frac{1}{2} I+ K' \right)\right)\sigma &= -\hat{\epsilon}\sum_i^{N_q} q_i \frac{\partial G}{\partial n}\label{eq:pcm}
\end{align}
where $G$ is the Laplace Green's function, $\hat{\epsilon} = (\epsilon_2 - \epsilon_1)/\epsilon_2$ and $K'$ is the normal electric field operator~\cite{Bardhan09_disc}.  Eq.~\ref{eq:pcm} is well known under any of several names, e.g. polarizable continuum model (PCM)~\cite{Miertus81} and apparent-surface charge (ASC) method~\cite{Cammi95}.  The surface charge $\sigma(r)$ induces a Coulomb potential called the \textit{reaction potential} $\varphi^{REAC}(r)$, because it arises from solvent polarization, and the quantity of interest, the solute-solvent interaction energy is $\frac{1}{2}\int \rho(r) \varphi^{REAC}(r) dr$.





\section{A Multiscale Model Incorporating Nonlocal Dielectric Behavior}\label{sec:nonlocal}


\begin{figure}[ht]
\centering
\includegraphics[width=0.6\textwidth]{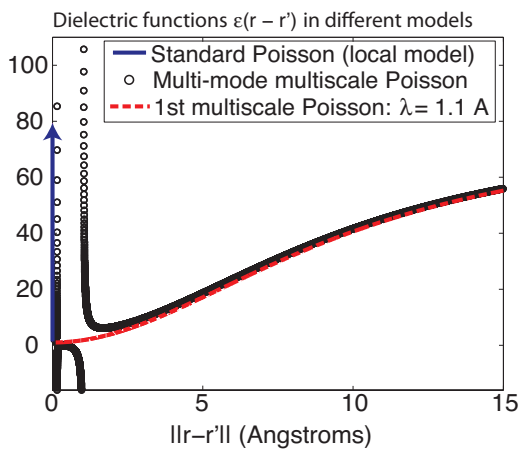}
  \caption{Multiscale continuum theory with nonlocal dielectric response.} \protect\label{fig:multiscale}
\end{figure}

Standard continuum models treat water as a macroscopic dielectric material; however, water molecules are not point particles or point dipoles. Macroscopic models miss important correlations induced by finite-size effects such as solvent molecules' finite
size~\cite{Roux99}, and hydrogen bonding and solvent structure, e.g. bridging solvent molecules.  Fully atomistic molecular-dynamics (MD) simulations correctly reproduce the main details of the length-scale dependence of water's dielectric response~\cite{Bopp98} (Fig.~\ref{fig:multiscale}), including charge oscillations, and we are seeking to implement these details using a computationally tractable BEM.

Macroscopic dielectric models derive from Gauss's law relating the electric flux $\mathbf{D}(\vr)$ to a fixed charge distribution $\rho(\vr)$, $\nabla \cdot \mathbf{D}(\vr) = \rho(\vr)/\epsilon_0$.
The Poisson equation follows by specifying the medium's relationship between the potential and the flux, traditionally $\mathbf{D}(\vr) = \epsilon_0 \epsilon(\vr)\mathbf{E}(\vr)$. Similar to gradient-elasticity theories~\cite{Bardhan13_nonlocal_review}, our multiscale theories improve on macrosopic models by creating a nonlocal relationship between the two: \mbox{$\mathbf{D}(\vr) = \int \epsilon_0 \epsilon(\vr-\vr') \mathbf{E}(\vr') d\vr'$}. The simplest version, the Lorentz nonlocal dielectric model~\cite{Hildebrandt04}, models dielectric correlations that decay with a characteristic length $\lambda_W$ from the short-range optical permittivity $\epsilon_{\infty}$,
\begin{align}
  \epsilon(\vr-\vr') & = \epsilon_\infty \delta(\vr - \vr') + \frac{\epsilon_{w} - \epsilon_\infty}{\lambda^2_W}
    \frac{e^{-\frac{|\vr - \vr'|}{\lambda_W}}}{4\pi|\vr - \vr'|}.\label{eq:lorentz-dielectric-function}
\end{align}
Because nonlocal models lead to integrodifferential equations of the form
\begin{align}
  \nabla \cdot \int \epsilon(\vr-\vr') \nabla \varphi(\vr') d\vr' = -\rho(\vr), 
\end{align}
which much harder to solve than standard Poisson PDE, the simple Lorentz nonlocal model is the only one applied to date in studies of polyatomic biomolecules~\cite{Bardhan11_pka}.  

Hildebrandt et al. made a key observation in noticing that the second term of Eq.~\ref{eq:lorentz-dielectric-function}, \mbox{$\exp(|\vr-\vr'|/\lambda)/|\vr-\vr'|$}, is the free-space Yukawa Green's function~\cite{Hildebrandt04}. Thus, the second component of $\mathbf{D}(\vr)$ solves a Yukawa equation, and applying the Helmholtz decomposition allows the nonlocal theory to be written as a coupled PDE system, after introducing an auxiliary potential $\psi_\mathrm{I}$,
\begin{align}
  \nabla^2 \varphi_\mathrm{IO}(\vr) &= -\rho(\vr), & \vr \in \mathrm{region\;II}\label{eq:coupled-nonlocal}\\
  \nabla^2 \psi_\mathrm{I}(\vr) &= 0, & \vr \in \mathrm{region\;I,II}\nonumber\\
  \left(\nabla^2 - \frac{1}{\Lambda^2}\right) \varphi_\mathrm{I}(\vr) &= -\frac{1}{\lambda^2} \psi_\mathrm{I}(\vr),  &\vr \in \mathrm{region\;I}\nonumber
\end{align}
where $\Lambda = \lambda_W\sqrt{\epsilon_\infty /   \epsilon_{\mathrm{w}}}$. The exact displacement boundary
condition also becomes nonlocal, but model studies show that this nonlocality can be safely omitted in many calculations as its impact is small~\cite{Weggler_thesis}, allowing use of the approximate, purely local boundary conditions~\cite{Hildebrandt04}
\begin{align}
  \varphi_\mathrm{II}(\vr) \mid_\Gamma &= \varphi_\mathrm{I}(\vr) \mid_\Gamma, \\
  \epsilon_0 \epsilon_\mathrm{p} \frac{\partial}{\partial n} \varphi_\mathrm{II}(\vr) \mid_\Gamma &= \frac{\partial}{\partial n} \psi_\mathrm{I}(\vr) \mid_\Gamma, \\
  \frac{\partial}{\partial n}\psi_\mathrm{II}(\vr) \mid_\Gamma &= \epsilon_0\epsilon_\infty \frac{\partial}{\partial n} \varphi_\mathrm{I}(\vr) \mid_\Gamma.
\end{align}
A change of variables $\Psi = \frac{1}{\epsilon_\infty}\left(\frac{1}{\epsilon_0} \psi_\mathrm{II} - \epsilon_\mathrm{p} \varphi_\mathrm{mol}\right)$ improves scaling, and repeated applications of Green's theorem lead to the BIE system
\begin{align}
  \begin{sbmatrix}{ccc}
    \frac{1}{2}-\KYL & -\frac{\epsilon_{\mathrm{p}}}{\epsilon_\infty}\VYL-\frac{\epsilon_{\mathrm{p}}}{\epsilon_{\mathrm{w}}}\VDR_\Lambda & \frac{\epsilon_\infty}{\epsilon_{\mathrm{w}}}\KDR_\Lambda\\
    \frac{1}{2}+\KL & -\VL & \\
                    & \frac{\epsilon_{\mathrm{p}}}{\epsilon_\infty} \VL & \frac{1}{2}-\KL
  \end{sbmatrix}
  \begin{sbmatrix}{c}
    \varphi_{\mathrm{II}} \\ \frac{\partial\varphi_{\mathrm{II}}}{\partial n} \\ \Psi
  \end{sbmatrix}
=
  \begin{sbmatrix}{c} \xi \\ 0 \\ 0 \end{sbmatrix},\label{eq:full-formulation}
\end{align}
where $V^{L,Y}$ and $K^{L,Y}$ are the single
and double layer operators for the Laplace and Yukawa kernels,
 $V^{DR}_\Lambda = V^Y_\Lambda - V^L$ and $K^{DR}_\Lambda = K^Y_\Lambda - K^L$, and
\begin{align}
  \xi = -\left(\frac{1}{2}-\KYL+\frac{\epsilon_{\mathrm{p}}}{\epsilon_{\mathrm{w}}}\KDR_\Lambda\right) \varphi_\mathrm{mol} - \left(\frac{\epsilon_{\mathrm{p}}}{\epsilon_\infty}\VYL - \frac{\epsilon_{\mathrm{p}}}{\epsilon_{\mathrm{w}}}\VDR_\Lambda\right)\frac{\partial\varphi_\mathrm{mol}}{\partial   n}.
\end{align}
For theory to match measurements of protein pH-dependence~\cite{Nielsen11}, standard models need to set $\epsilon_{\mathrm{p}} > 10$, well beyond experiments indicating that $\epsilon_{\mathrm{p}} < 6$, e.g.~\cite{Kukic13}. We have proposed that nonlocal solvent response decreases dielectric contrast similarly to increasing $\epsilon_p$, allowing use of experimentally reasonable $\epsilon_p$~\cite{Bardhan11_pka}.  For example, for a spherical boundary, all of the boundary-integral operators are
diagonalized by the surface spherical harmonics~\cite{Xie12,Bardhan14_analytical_nonlocal}; in~\cite{Bardhan14_analytical_nonlocal} we used these
expansions to study nonlocal response for realistic protein charge distributions rapidly and with high accuracy.
\begin{figure}[ht]
\centering    \includegraphics[width=\textwidth]{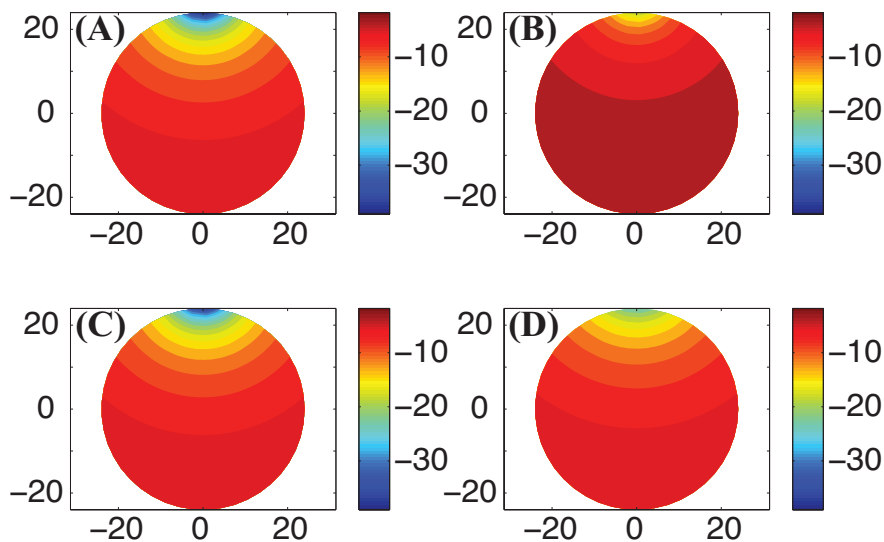}
\caption{Comparison of local and nonlocal excitations on the sphere.}
\protect  \label{fig:reac-pot-vary-epsilon-surface-charge}
\end{figure}
Fig.~\ref{fig:reac-pot-vary-epsilon-surface-charge} (from~\cite{Bardhan14_analytical_nonlocal}) illustrates key differences between macroscopic local-response theory and nonlocal theory for a spherical protein (radius~24~\AA) with a charge 2~\AA~from the solvent interface. All potentials are plotted on the same scale and $\epsilon_\mathrm{w} = 80$ for all models. The panels show (a) local model, $\epsilon_\mathrm{p} = 2$;
(b) local model, $\epsilon_\mathrm{p} = 4$; (c) nonlocal model,
$\epsilon_\mathrm{p} = 2$ and $\lambda = 1$~\AA; (d) nonlocal model, $\epsilon_\mathrm{p} = 2$ and $\lambda = 10$~\AA.  In addition to analytical methods, we have implemented a fast BEM solver for these
equations~\cite{Bardhan11_DAC}.

Our calculations of amino acid titrations shows that nonlocality reduces dielectric contrast, which may be of functional importance~\cite{Bardhan11_pka}.  However, the Lorentz model fails to reproduce obvious features of solvent behavior~\cite{Attard90}, such as the presence of charge-density oscillations (Fig.~\ref{fig:multiscale}).  To test such a model, we have derived a similar coupled-PDE decomposition for solving the charge-oscillation nonlocal model of Fig.~\ref{fig:multiscale}~\cite{Bardhan13_nonlocal_review}.  Current work focuses on deriving a BIE formulation to solve it efficiently. Many important open questions remain in nonlocal theory, including uniqueness and provable bounds on the errors of approximating the boundary condition~\cite{Bardhan13_nonlocal_review}.




\section{A Nonlinear Boundary Condition for Atomistic  Effects at Biomolecule Surfaces}\label{sec:nlbc}

\begin{figure}[ht]
    \includegraphics[width=\textwidth]{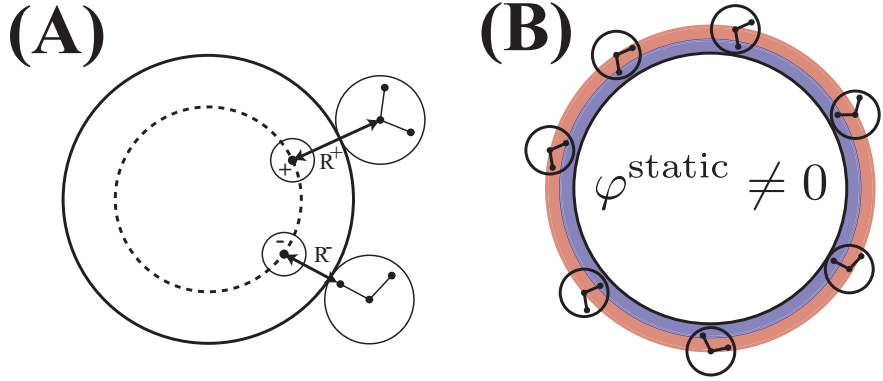}
    \caption{Charge-hydration asymmetry results from distinct mechanisms. (A) Negative surface charges interact more strongly because smaller water hydrogens can approach more closely. (B) Water-water correlations around the uncharged molecule create a     nearly uniform dipole charge layer.}\protect\label{fig:asymmetry-mechanisms}
\end{figure}

\begin{figure}[h]
  \includegraphics[width=\textwidth]{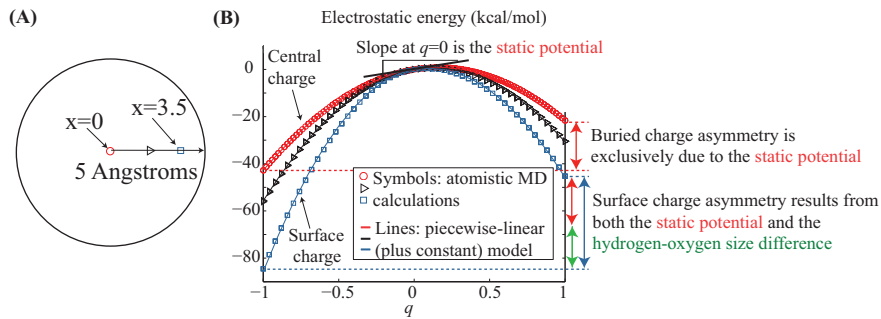}
  \caption{Atomistic simulations demonstrate that distinct mechanisms of asymmetry affect buried and surface charges to different degrees (adapted from~\cite{Bardhan12_asymmetry}). (A) An approximately spherical solute of radius 5 \AA. Molecular dynamics simulations were performed for a single charge at the center or closer to the surface. (B) Results from MD simulations (symbols) and a piecewise-linear model (curves).  At $q=0$, the slopes for all of the curves are equal to the static (surface) potential.}
  \label{fig:asymmetry-sphere-fits}
\end{figure}

As noted previously, water's finite size leads to correlations, and these become more pronounced at surfaces. As shown in Fig.~\ref{fig:asymmetry-mechanisms}, the smaller hydrogens can approach solute charges more closely than the larger oxygen~(from~\cite{Bardhan12_asymmetry}); as a result, negative charges interact more strongly with solvent~(Fig.~\ref{fig:asymmetry-sphere-fits}).  Furthermore, correlations between waters at the surface produce a dipole charge layer, creating a large, nearly constant potential $\phi^{static}$ in the interior even when the solute is devoid of charges.  As a result, a charge's energy is better modeled as $\frac{1}{2} L_+ q^2 + \phi^{static} q$ for $q > 0$, and as $\frac{1}{2} L_- q^2 + \phi^{static} q$ for $q < 0$, where $L_- < L_+$ (both are negative).
However, standard Poisson models give energies that are symmetric with respect to the charge's sign, and can only account for charge-hydration asymmetry by detailed parameterization, i.e. changing dozens of atomic radii until calculations fit a set of reference energies.  The resulting radii often run counter to physical intuition. For example, in one parameter set~\cite{Nina97}, the radii of carbon atoms range from 2.04~\AA~to 2.86~\AA, depending on its chemical bonds. However, changing radii does not directly address physics.

\begin{figure}[h]
  \includegraphics[width=\textwidth]{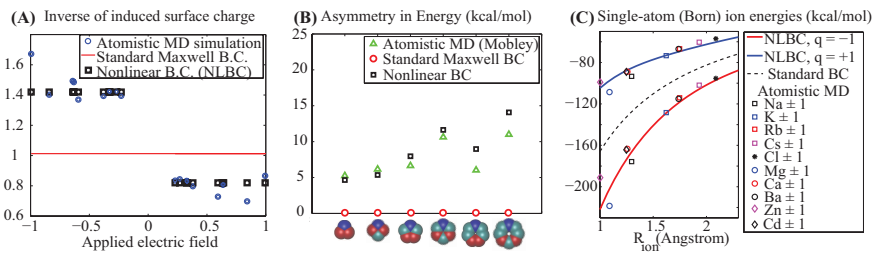}
  \caption{Comparison with atomistic molecular-dynamics (MD) simulations validates  our NLBC formulation for single atoms and challenging model problems. (A), The Standard Maxwell Boundary Condition (SMBC) fails to reproduce the surface charge obtained from MD, whereas our nonlinear boundary condition (NLBC) produces the correct qualitative picture.  (B) Our NLBC gives semi-quantitative  agreement with the MD simulations of Mobley et al. for asymmetry in model problems; the SMBC gives 0 exactly. (C) The NLBC is  quantitatively accurate for atomic ions \textit{without} fitting their atomic radii (modified from~\cite{Bardhan14_asym,Bardhan15_PIERS}).}
  \label{fig:nlbc-motivation}
\end{figure}

In contrast, we recently modeled asymmetry by adding a nonlinear correction to the standard Maxwell boundary condition (SMBC)~\cite{Bardhan14_asym}. Our modified boundary condition directly changes response depending on the local electric field at the surface, mirroring the actual physics.  Consider that the standard Maxwell boundary condition $\epsilon_w \frac{\partial \varphi_{I}}{\partial n}=\epsilon_p \frac{\partial \varphi_{II}}{\partial n}$ leads to a surface charge satisfying
\begin{align}
  \frac{\sigma(\vr_\Gamma)}{\epsilon_1} = \frac{\partial\varphi_1}{\partial n}(\vr_\Gamma) - \frac{\partial\varphi_2}{\partial n}(\vr_\Gamma),\label{eq:ecf-bc}
\end{align}
which explicitly reveals the sign-symmetry of the induced charge. On the other hand, a surface charge matching MD results indicates nonlinear response at the boundary, namely that a negative charge provides slightly greater surface charge (Fig.~\ref{fig:nlbc-motivation}(A)).  We model this sign dependence with the phenomenological nonlinear boundary condition (NLBC)
\begin{align}
  \left(1 + f(e_n)\right) \frac{\sigma(\vr_\Gamma)}{\epsilon_1} = \frac{\partial\varphi_1}{\partial n}(\vr_\Gamma) - \frac{\partial\varphi_2}{\partial n}(\vr_\Gamma),\label{eq:ecf-nlbc}
\end{align}
where $E_n$ is the normal electric field just \textit{inside} the boundary $\Gamma$ and $f$ is a smoothed step function,
\begin{align}
  f(E_n) &= \frac{\epsilon_1}{\epsilon_2-\epsilon_1} - h(E_n),\\
  h(E_n) &= \alpha \tanh(\beta E_n - \gamma) + \mu.\label{eq:tanh}
\end{align}
Eliminating $\sigma$ gives 
\begin{align}
  \frac{f}{1+f} \frac{\partial \varphi_1}{\partial n}(\vr_\Gamma) = \frac{\partial \varphi_2}{\partial n}(\vr_\Gamma).\label{eq:nlbc}
\end{align}
The parameters have physical meaning: $\alpha$ models the magnitude of asymmetry, $1/\beta$ models the electric
field strength for saturation of the NLBC, and $\gamma$ models water's intrinsic preference for one orientation over
another.  These parameters were successfully fit to reproduce dozens of MD free-energy calculations of Mobley et.
al.~\cite{Mobley08} on net-neutral fictitious ``bracelet'' and ``rod'' molecules~\cite{Bardhan14_asym} (one test set is shown in Figure~\ref{fig:nlbc-motivation}(B)).  For a pure water solvent, the NLBC Eq.~\ref{eq:nlbc} leads to a modified form of the PCM/ASC Eq.~\ref{eq:pcm},
\begin{align}
  \left(I  + \hat{\epsilon} \left(-\frac{1}{2} I+ K' \right) + h(E_n) \right)\sigma &= -\hat{\epsilon}\sum_i^{N_q} q_i \frac{\partial G}{\partial n}
\end{align}
where $E_n = -\sum_i q_i \frac{\partial G}{\partial n} - K \sigma$ is the interior field at the boundary.
Numerous previous nonlinear Poisson models fail to address charge-hydration asymmetry because they focus on saturation at high field strengths/charge densities, e.g. the nonlinear Poisson--Boltzmann equation~\cite{Sharp90} and dielectric saturation~\cite{Alper90,Hu12_Wei_nonlinear_Poisson_BJ}. However,
charge-hydration asymmetry contributes significant energetics even for low charge densities and neutral molecules~\cite{Bardhan12_asymmetry,Bardhan14_asym}. Moreover, most nonlinear models are still charge-sign symmetric (though not all~\cite{Hu12_Wei_nonlinear_Poisson_BJ}).

\begin{figure}[h]
  \includegraphics[width=\textwidth]{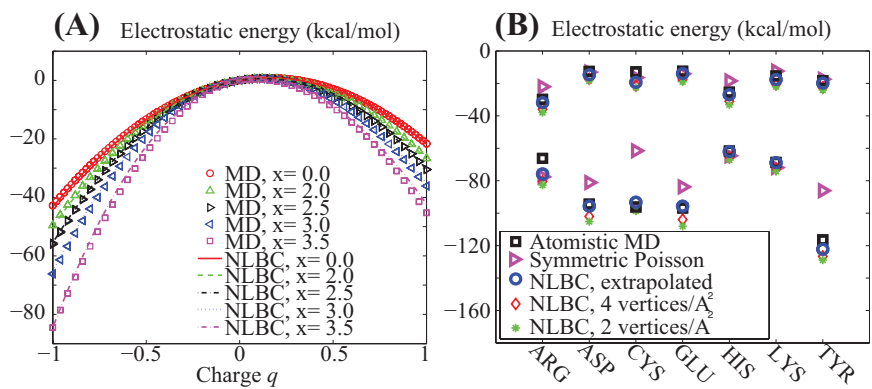}
  \caption{Comparison with explicit MD simulations validates the accuracy of our NLBC model for both asymmetry mechanisms and for complex biomolecules. (A) NLBC results accurately reproduce asymmetry for both buried and surface charges, see Fig.~\ref{fig:asymmetry-sphere-fits}(B). (B) Energies for protonated and deprotonated forms of titratable amino (modified from~\cite{Bardhan14_asym,Bardhan15_PIERS}): ARG=arginine, ASP=aspartic acid, CYS=cysteine, GLU=glutamic acid, HIS=histidine, LYS=lysine, TYR=tyrosine.}
  \label{fig:nlbc-results}
\end{figure}

Our asymmetric Poisson model reproduces explicit-solvent calculations for Born ions, regardless of their natural charge (Figure~\ref{fig:nlbc-motivation}(C))~\cite{Bardhan14_asym}. We emphasize that unlike most models, we have not fit the ion radii. Our radii are simply the MD radii, scaled by the scale factor 0.9 (the same scale factor obtained by parameterization
against the Mobley simulations~\cite{Mobley08}).  Our NLBC also correctly predicts charging free energies for the model problem in Figure~\ref{fig:nlbc-results}(A), a
charge embedded in a 5-\AA-radius ~sphere.  This result is particularly important because our model correctly predicts the
crossover between the interface-potential asymmetry (for buried charges) and the NLBC asymmetry (for surface charges).
Calculations of the protonated and deprotonated forms of the titratable amino acids illustrate that the NLBC model works well for atomistic models of more complicated molecules containing polar and charged chemical groups~\cite{Bardhan14_asym} (Fig.~\ref{fig:nlbc-results}(B)). The results in the Figure indicate that our model can reproduce more expensive MD simulations \textit{better than a competing symmetric-response theory with dozens   more fitting parameters}. Moreover, because the NLBC parameters have a physical basis, it may soon be possible to calculate them independently from first principles.





\section{Conclusion}\label{sec:discussion}

Electrostatic interactions between molecules and surrounding fluid represent a challenging problem in biology and chemistry~\cite{Sharp90}. As a consequence of the success of continuum models for these interactions, some of the most cited BEM papers address the polarizable continuum model (PCM) for the problem~\cite{Miertus81}.  Emerging areas for biomolecular modeling are challenging popular PDE solvers and heuristic electrostatic models, demonstrating the need for more realistic theories and fast, accurate simulations and opening the door for advanced BIE approaches. For reasons of space, we have not been able to discuss many other advances, such as interactions between thousands of proteins~\cite{Yokota11} or a scalable approximation theory based on low-rank approximation of the integral operators~\cite{Bardhan08_BIBEE,Bardhan09_bounds,Bardhan11_Knepley}. In this paper, we have focused a few examples of the modeling opportunities and impact that BIE theory can have in this area of computational science, and we hope that the BEM community will contribute its expertise to the many open questions that remain.

\bibliography{implicit-review}
\bibliographystyle{witpress}

\end{document}

%% file: preamble.tex
\usepackage{graphicx}
\usepackage{amssymb}
\usepackage{amsmath}
\usepackage{color}
\usepackage{url}
\usepackage{multirow}
\usepackage{graphicx} 
\usepackage{wrapfig,floatrow} 

\newenvironment{sbmatrix}[1]{\left[ \begin{array}{#1}}%
	{\end{array} \right]}
\def\KL{K^L}
\def\VL{V^L}
\def\KYL{K^L_\Lambda}
\def\VYL{V^L_\Lambda}
\def\KDR{K^{DR}}
\def\VDR{V^{DR}}

\def\vr{\mathbf{r}}

\usepackage[colorinlistoftodos, textwidth=4cm, shadow]{todonotes}

\newcommand{\matt}[1]{\todo[color=green!40,inline]{{\bf Matt:} #1}}

\newcommand{\arxiv}{\let\arxivFig\matt}

\def\XXint#1#2#3{{\setbox0=\hbox{$#1{#2#3}{\int}$}
     \vcenter{\hbox{$#2#3$}}\kern-.5\wd0}}

%% file: bem38.bbl
\begin{thebibliography}{10}

\bibitem{Bardhan12_review}
Bardhan, J.P., Biomolecular electrostatics---{I} want your solvation (model).
  \emph{Computational Science and Discovery}, \textbf{5}, p. 013001, 2012.

\bibitem{Roux99}
Roux, B. \& Simonson, T., Implicit solvent models. \emph{Biophys Chem},
  \textbf{78}, pp. 1--20, 1999.

\bibitem{Miertus81}
Miertus, S., Scrocco, E. \& Tomasi, J., Electrostatic interactions of a solute
  with a continuum -- a direct utilization of \textit{ab initio} molecular
  potentials for the prevision of solvent effects. \emph{Chem Phys},
  \textbf{55(1)}, pp. 117--129, 1981.

\bibitem{Bardhan13_nonlocal_review}
Bardhan, J.P., Gradient models in molecular biophysics: progress, challenges,
  opportunities. \emph{Journal of Mechanical Behavior of Materials},
  \textbf{22}, pp. 169--184, 2013.

\bibitem{Bardhan15_CEBA}
Bardhan, J.P., Boundary-integral and boundary-element methods for biomolecular
  electrostatics: progress, challenges, and important lessons from ceba 2013.
  \emph{Computational Electrostatics for Biological Applications}, Springer,
  2015.

\bibitem{Sharp90}
Sharp, K.A. \& Honig, B., Electrostatic interactions in macromolecules: Theory
  and applications. \emph{Annu Rev Biophys Bio}, \textbf{19}, pp. 301--332,
  1990.

\bibitem{Bardhan09_disc}
Bardhan, J.P., Numerical solution of boundary-integral equations for molecular
  electrostatics. \emph{J Chem Phys}, \textbf{130}, p. 094102, 2009.

\bibitem{Cammi95}
Cammi, R. \& Tomasi, J., Remarks on the apparent-surface charges ({ASC})
  methods in solvation problems: {Iterative} versus matrix-inversion procedures
  and the renormalization of the apparent charges. \emph{J Comput Chem},
  \textbf{16(12)}, pp. 1449--1458, 1995.

\bibitem{Bopp98}
Bopp, P.A., Kornyshev, A.A. \& Sutmann, G. \emph{J Chem Phys}, \textbf{109}, p.
  1939, 1998.

\bibitem{Hildebrandt04}
Hildebrandt, A., Blossey, R., Rjasanow, S., Kohlbacher, O. \& Lenhof, H.P.,
  Novel formulation of nonlocal electrostatics. \emph{Phys Rev Lett},
  \textbf{93}, p. 108104, 2004.

\bibitem{Bardhan11_pka}
Bardhan, J.P., Nonlocal continuum electrostatic theory predicts surprisingly
  small energetic penalties for charge burial in proteins. \emph{J Chem Phys},
  \textbf{135}, p. 104113, 2011.

\bibitem{Weggler_thesis}
Weggler, S., \emph{Correlation induced electrostatic effects in biomolecular
  systems}. Ph.D. thesis, Universit\"{a}t des Saarlandes, 2010.

\bibitem{Nielsen11}
Nielsen, J.E., Gunner, M., Garc{\'\i}a-Moreno, E. et~al., The pka cooperative:
  A collaborative effort to advance structure-based calculations of pka values
  and electrostatic effects in proteins. \emph{Proteins: Structure, Function,
  and Bioinformatics}, \textbf{79(12)}, pp. 3249--3259, 2011.

\bibitem{Kukic13}
Kukic, P., Farrell, D., {McIntosh}, L.P., E., B.G.M., Jensen, K.S., Toleikis,
  Z., Teilum, K. \& Nielsen, J.E., Protein dielectric constants determined from
  {NMR} chemical shift parameters. \emph{J Am Chem Soc}, 2013.

\bibitem{Xie12}
Xie, D., Jiang, Y., Brune, P. \& Scott, L.R., A fast solver for a nonlocal
  dielectric continuum model. \emph{SIAM Journal of Scientific Computing},
  \textbf{34}, pp. B107--B126, 2012.

\bibitem{Bardhan14_analytical_nonlocal}
Bardhan, J.P., Brune, P.R. \& Knepley, M.G., Analytical nonlocal electrostatics
  using eigenfunction expansion of boundary-integral operators. \emph{Molecular
  Based Mathematical Biology}, \textbf{3}, 2015.

\bibitem{Bardhan11_DAC}
Bardhan, J.P. \& Hildebrandt, A., A fast solver for nonlocal electrostatic
  theory in biomolecular science and engineering. \emph{IEEE/ACM Design
  Automation Conference (DAC)}, 2011.

\bibitem{Attard90}
Attard, P., Wei, D. \& Patey, G.N., Critical comments on the nonlocal
  dielectric function employed in recent theories of the hydration force.
  \emph{Chemical Physics Letters}, \textbf{172}, pp. 69--72, 1990.

\bibitem{Bardhan12_asymmetry}
Bardhan, J.P., Jungwirth, P. \& Makowski, L., Affine-response model of
  molecular solvation of ions: Accurate predictions of asymmetric charging free
  energies. \emph{J Chem Phys}, \textbf{137}, p. 124101, 2012.

\bibitem{Nina97}
Nina, M., Beglov, D. \& Roux, B., Atomic radii for continuum electrostatics
  calculations based on molecular dynamics free energy simulations. \emph{J
  Phys Chem B}, \textbf{101}, pp. 5239--5248, 1997.

\bibitem{Bardhan14_asym}
Bardhan, J.P. \& Knepley, M.G., Modeling charge-sign asymmetric solvation free
  energies with nonlinear boundary conditions. \emph{J Chem Phys},
  \textbf{141}, p. 131103, 2014.

\bibitem{Mobley08}
Mobley, D.L., Dill, K.A. \& Chodera, J.D., Treating entropy and conformational
  changes in implicit solvent simulations of small molecules. \emph{J Phys Chem
  B}, \textbf{112}, pp. 938--946, 2008.

\bibitem{Alper90}
Alper, H.E. \& Levy, R.M., Field strength dependence of dielectric saturation
  in liquid water. \emph{J Phys Chem}, \textbf{94}, pp. 8401--8403, 1990.

\bibitem{Hu12_Wei_nonlinear_Poisson_BJ}
Hu, L. \& Wei, G.W., Nonlinear {Poisson} equation for heterogeneous media.
  \emph{Biophys J}, \textbf{103}, pp. 758--766, 2012.

\bibitem{Bardhan15_PIERS}
Bardhan, J.P., Tejani, D., Wieckowski, N., Ramaswamy, A. \& Knepley, M.G., A
  nonlinear boundary condition for continuum models of biomolecular
  electrostatics. \emph{Progress in Electromagnetics Research Symposium
  (PIERS)}, 2015.

\bibitem{Yokota11}
Yokota, R., Bardhan, J.P., Knepley, M.G., Barba, L.A. \& Hamada, T.,
  Biomolecular electrostatics using a fast multipole {BEM} on up to 512 {GPUs}
  and a billion unknowns. \emph{Comput Phys Commun}, \textbf{182}, pp.
  1272--1283, 2011.

\bibitem{Bardhan08_BIBEE}
Bardhan, J.P., Interpreting the {Coulomb}-field approximation for
  {Generalized}-{Born} electrostatics using boundary-integral equation theory.
  \emph{J Chem Phys}, \textbf{129(144105)}, 2008.

\bibitem{Bardhan09_bounds}
Bardhan, J.P., Knepley, M.G. \& Anitescu, M., Bounding the electrostatic free
  energies associated with linear continuum models of molecular solvation.
  \emph{J Chem Phys}, \textbf{130}, p. 104108, 2009.

\bibitem{Bardhan11_Knepley}
Bardhan, J.P. \& Knepley, M.G., Mathematical analysis of the boundary-integral
  based electrostatics estimation approximation for molecular solvation: Exact
  results for spherical inclusions. \emph{J Chem Phys}, \textbf{135}, p.
  124107, 2011.

\end{thebibliography}
